\documentclass{elsarticle}
\usepackage{amsfonts}

\usepackage{graphicx}
\usepackage{pstricks}
\usepackage{amsmath}
\usepackage{amsxtra}
\def\dim{\operatorname{dim}}

\newtheorem{theorem}{Theorem}[section]
\newtheorem{lemma}[theorem]{Lemma}
\newproof{proof}[theorem]{Proof}
\newtheorem{proposition}[theorem]{Proposition}
\newtheorem{corollary}[theorem]{Corollary}

\newtheorem{example}[theorem]{Example}
\newtheorem{remark}[theorem]{Remark}

\numberwithin{equation}{section}
\begin{document}

%\linenumbers

\begin{frontmatter}

\title{Cubic column relations in truncated moment problems}

\author{Ra\'{u}l E. Curto \footnote{The first named author was partially supported by NSF Grants
DMS-0400741 and DMS-0801168.}}
\address{Department of Mathematics, The University of Iowa, Iowa City, Iowa
52242}
\ead{raul-curto@uiowa.edu}
\ead[url]{http://www.math.uiowa.edu/\symbol{126}rcurto/}

\author{Seonguk Yoo \footnote{The second named author was partially supported by a University of Iowa Graduate
College Summer Fellowship and by a 2007 IMA PI Summer Program Graduate Student Fellowship (held at Texas A \& M University).}}
\address{Department of Mathematics, The University of Iowa, Iowa City, Iowa 52242}
\ead{seonguk-yoo@uiowa.edu}

\begin{abstract}

For the truncated moment problem associated to a complex sequence $\gamma
^{(2n)}=\{\gamma _{ij}\}_{i,j\in Z_{+},i+j \leq 2n}$ to have a representing measure $\mu $, it is necessary for the moment matrix $M(n)$ to be positive semidefinite, and for the
algebraic variety $\mathcal{V}_{\gamma}$ to satisfy $\operatorname{rank}\;M(n) \leq \;$ card$\;\mathcal{V}_{\gamma}$ as
well as a consistency condition: the Riesz functional vanishes on every polynomial of degree at most $2n$ that vanishes on $\mathcal{V}_{\gamma}$. \ In previous work
with L. Fialkow and M. M\"{o}ller, the first-named author proved that for the extremal case (rank$\;M(n)=$ card$\;\mathcal{V}_{\gamma}$), positivity and consistency are sufficient for the existence of a
representing measure. \ 

In this paper we solve the truncated moment problem for \textit{cubic} column relations in $M(3)$ of the form $Z^{3}=itZ+u\bar{Z}$ ($u,t \in \mathbb{R}$); we do this by checking consistency. \ For $(u,t)$ in the open cone determined by $0 < \left|u\right| < t < 2 \left|u\right|$, we first prove that the algebraic variety has exactly $7$ points and $\operatorname{rank}\;M(3)=7$; we then apply the above mentioned result to obtain a concrete, computable, necessary and sufficient condition for the existence of a representing measure. 
\end{abstract}


\end{frontmatter}

\section{\label{Introduction}Introduction}

Given a collection of complex numbers $\gamma \equiv \gamma ^{(2n)}:\gamma
_{00},\gamma _{01},\gamma _{10},\cdots ,$ $\gamma _{0,2n},$\ $\gamma
_{1,2n-1},\cdots ,\gamma _{2n-1,1},\gamma _{2n,0},$ with $\gamma _{00}>0$
and $\gamma _{ji}=\bar{\gamma}_{ij},$ the \textit{truncated complex
moment problem} (TCMP) consists of finding a positive Borel measure $\mu $
supported in the complex plane $\mathbb{C}$ such that $\gamma _{ij}=\int 
\bar{z}^{i}z^{j}\;d\mu \;\;\;(0\leq i+j\leq 2n);$ $\gamma $ is called a 
\textit{truncated moment sequence} (of order $2n$) and $\mu $ is called a 
\textit{representing measure} for $\gamma $. \ If, in addition, we require $%
\operatorname{supp}\;\mu \subseteq K\;$(closed) $\subseteq \mathbb{C}$, we speak of
the $K$-TCMP. $\ $Naturally associated with each TCMP is a moment matrix $%
M(n)\equiv M(n)(\gamma )$, whose concrete construction will be given in Subsection %
\ref{TCMP}. $\ M(n)$ detects the positivity of the \textit{Riesz functional} $%
p \mapsto \sum_{ij}a_{ij}\gamma _{ij}\;\;(p(z,\bar{z})\equiv
\sum_{ij}a_{ij}\bar{z}^{i}z^{j})$ on the cone generated by the collection $%
\{p\bar{p}:p\in \mathbb{C}[z,\bar{z}]\}$. $\ $In addition to its importance
for applications, a complete solution of TCMP would readily lead to a
solution of the \textit{full} moment problem, via a weak-* convergence argument, as
shown by J. Stochel \cite{Sto2}. \ While we primarily focus on truncated moment problems, the full moment problem (in one or several variables) has been widely studied; see, for example, \cite{BaTe}, \cite{blekherman}, \cite{Dem}, \cite{KuMa}, \cite{Lau2}, \cite{Lau3}, \cite{PoSc}, \cite{Pu1}, \cite{PuSch}, \cite{PuSchm}, \cite{PuVa1}, \cite{PuVa2}, \cite{Rie}, \cite{Sche1}, \cite{Sche2}, \cite{Sch1}, \cite{Sch2}, \cite{Sch4}, \cite{StSz1}, \cite{polar}, \cite{Vas2}.

Building on previous work for the case of \textit{real} moments, several years ago the first named author and L. Fialkow introduced in 
\cite{tcmp1}, \cite{tcmp2} and \cite{tcmp3} an approach to TCMP based on
matrix positivity and extension, combined with a new ``functional calculus''
for the columns of $M(n)$ (labeled $1,Z,\bar{Z},Z^{2},\bar{Z}Z,\bar{Z}^{2},$%
...). \ This allowed them to show that TCMP is soluble in the following cases:

\noindent (i) TCMP is of \textit{flat data} type \cite{tcmp1}, i.e., $\operatorname{rank}\;M(n)=\operatorname{rank}\;M(n-1)$ (this case
subsumes all previous results for the Hamburger, Stieltjes, Hausdorff, and
Toeplitz truncated moment problems \cite{Houston});

\noindent(ii) the column $\bar{Z}$ is a linear combination of the columns $1$
and $Z$ \cite[Theorem 2.1]{tcmp2};

\noindent(iii) for some $k\leq \lbrack n/2]+1,$ the analytic column $Z^{k}$ is
a linear combination of columns corresponding to monomials of lower degree 
\cite[Theorem 3.1]{tcmp2};

\noindent(iv) the analytic columns of $M(n)$ are linearly dependent and span $%
\mathcal{C}_{M(n)}$, the column space of $M(n)$ \cite[Corollary 5.15]{tcmp1};

\noindent(v) $M(n)$ is singular and subordinate to conics \cite{tcmp5}, \cite%
{tcmp6}, \cite{tcmp7}, \cite{tcmp9};

\noindent(vi) $M(n)$ admits a rank-preserving moment matrix extension $M(n+1)$, i.e., an extension $M(n+1)$ which is flat 
\cite{tcmp10};

\noindent(vii) $M(n)$ is extremal, i.e., $\operatorname{rank}\;M(n)=\operatorname{card} \;%
\mathcal{V}(\gamma ^{(2n)})$, where $\mathcal{V}(\gamma )\equiv \mathcal{V}%
(\gamma ^{(2n)})$ is the algebraic variety of $\gamma $ (see Subsection \ref%
{extremal} below) \cite{tcmp11}.

The common feature of the above mentioned cases is
the presence, at the level of $\mathcal{C}_{M(n)}$, of algebraic conditions
implied by the existence of a representing measure with support in a proper
real algebraic subset of the plane. \ However, the chief attraction
of the truncated moment problem (TMP) is its naturalness: since the data set is
finite, we can apply ``finite'' techniques, grounded in finite dimensional
operator theory, linear algebra, and algebraic geometry, to develop
algorithms for explicitly computing finitely atomic representing measures. \ 

\subsection{\textbf{The Truncated Complex Moment Problem}} \label{TCMP}

Given a collection of complex numbers $\gamma \equiv \gamma ^{(2n)}:\gamma
_{00},\gamma _{01},\gamma _{10},\cdots ,$ $\gamma _{0,2n},$\ $\gamma
_{1,2n-1},\cdots ,\gamma _{2n-1,1},\gamma _{2n,0},$ with $\gamma _{00}>0$
and $\gamma _{ji}=\bar{\gamma}_{ij},$ the associated moment matrix $%
M(n)\equiv M(n)(\gamma )$ is built as follows.  
\begin{equation*}
M(n):=\left( 
\begin{array}{cccc}
M[0,0] & M[0,1] & \cdots & M[0,n] \\ 
M[1,0] & M[1,1] & \cdots & M[1,n] \\ 
\cdots & \cdots & \cdots & \cdots \\ 
M[n,0] & M[n,1] & \cdots & M[n,n]%
\end{array}%
\right),
\end{equation*}
where 
$$
M[i,j]:=\left( 
\begin{array}{cccc}
\gamma _{i,j} & \gamma _{i+1,j-1} & \cdots & \gamma _{i+j,0} \\ 
\gamma _{i-1,j+1} & \gamma _{i,j} & \cdots & \gamma _{i+j-1,1} \\ 
\vdots & \vdots & \ddots & \vdots \\ 
\gamma _{0,i+j} & \gamma _{1,i+j-1} & \cdots & \gamma _{j,i}%
\end{array}%
\right).
$$
Observe that each rectangular block $M[i,j]$ is Toeplitz (that is, constant along diagonals), and that $M(n+1)=\left(\begin{array}{cccc} M(n) & B \\ B^{*} & C \end{array} \right)$, for some matrices $B$ and $C$.  Moreover, the results in \cite{tcmp1} and \cite{tcmp3} imply that each \textit{soluble} TMP with finite algebraic variety is eventually extremal (see also \cite{finitevariety}).

It is well known that the positivity of the moment matrix $M(n)$ is a necessary condition for the existence of a representing measure \cite{tcmp1}; thus, we always assume $M(n) \ge 0$. \ To check the positivity of a prospective moment matrix $M(n+1)$ given the positivity of
$M(n)$, we need the following classical result.

\begin{theorem}
\label{smul}(Smul'jan's Theorem \cite{Smu}) Let $A,B,C$ be matrices of complex numbers, with $A$ and $C$ square matrices. \ Then  
\begin{equation*}
\left( 
\begin{array}{cc}
A & B \\ 
B^{\ast } & C%
\end{array}%
\right) \geq 0 \iff \left\{ 
\begin{array}{c}
A\geq 0 \\ 
B=AW \; (\textrm{for some~} W)\\ 
C\geq W^{\ast }AW%
\end{array}%
\right. .
\end{equation*}%
Moreover, $\operatorname{rank}\;$$\left( 
\begin{array}{cc}
A & B \\ 
B^{\ast } & C%
\end{array}%
\right) =$$\operatorname{rank}\;$$A \iff C=W^{\ast }AW.$ 
\end{theorem}

\begin{remark}
When the rank equality occurs in Theorem \ref{smul}, we say that \newline $
\left(
\begin{array}{cc}
A & B \\ 
B^{\ast } & C%
\end{array}
\right)
\equiv 
\left(
\begin{array}{cc}
A & AW \\ 
W^{\ast }A & W^{\ast }AW%
\end{array}
\right)
$
is a flat extension of $A$. \ Note that while flat extensions are in principle easy to build, given a moment matrix $A \equiv M(n)$ the block $W^{\ast }AW \equiv W^{\ast }M(n)W$ may or may not satisfy the structural property of being Toeplitz. \ This is precisely the difficulty posed in generating flat extensions of positive moment matrices.
\end{remark}

\subsection{\label{extremal}\textbf{Extremal Moment Problems}}

Given a polynomial $p(z,\bar{z}) \equiv \sum_{ij}a_{ij}\bar{z}^{i}z^{j}$ we let $p(Z,\bar{Z}):=\sum_{ij}a_{ij}\bar{Z}^{i}Z^{j}$ (so that $p(Z,\bar{Z}%
)\in \mathcal{C}_{M(n)}$), and we let $\mathcal{Z}(p)$
denote the zero set of $p$. \ The assignment $p \mapsto p(Z,\bar{Z})$ is what we call the ``functional calculus." \ We define the \textit{algebraic variety} of $\gamma $ by 
\begin{equation} \label{variety}
\mathcal{V} \equiv \mathcal{V}(\gamma ):=\bigcap {}_{p(Z,\bar{Z})=0,\; \deg \;p\leq n}\mathcal{Z}%
(p).
\end{equation}%
Observe that $p(Z,\bar{Z})=M(n)\widehat{p}$ (where 
$\widehat{p}$ denotes the vector of coefficients of $p$), so that 
$p(Z,\bar{Z})=0$ if and only if $\widehat{p} \in \operatorname{ker} M(n)$. 
\ If $\gamma $ admits a representing measure $\mu $,
then the rank, $r$, of the moment matrix $M(n)$ is always bounded above by the 
cardinality, $v$, of $\mathcal{V}(\gamma )$; one actually has $\operatorname{supp} \mu \subseteq \mathcal{V}(\gamma )$ and 
$r\leq \operatorname{card}\;\operatorname{supp}\;\mu \leq v$ \cite{tcmp3}. \ Further, in
this case, if $p$ is any polynomial of degree at most $2n$ such that $p|_{%
\mathcal{V}}\equiv 0$, then the Riesz functional $\Lambda $ satisfies $%
\Lambda (p)=\int p~d\mu =0$. \ In summary, the following three conditions are
clearly necessary for the existence of a representing measure for $\gamma
^{(2n)}$: 
\begin{eqnarray}
\text{(Positivity)}~~M(n)\geq 0 \label{C1} \\  
\text{(Consistency)}~~p\in \mathcal{P}_{2n},\;p|_{\mathcal{V}}\equiv
0\Longrightarrow \Lambda (p)=0  \label{C2}\\ 
\text{(Variety Condition)}~~r\leq v\text{, i.e., }\operatorname{rank}\;M(n)\leq 
\operatorname{card}\;\mathcal{V}\text{.}  \label{C3}
\end{eqnarray}
The main result in \cite{tcmp11} shows that these three conditions are
indeed sufficient in the \textit{extremal} case ($r=v$). \ It was also proved in 
\cite{tcmp11} that Consistency cannot be replaced by the weaker condition that 
$M(n)$ is recursively generated, even if the algebraic variety is a planar cubic. \ ($M(n)$ is \textit{recursively generated} if for any column relation in $M(n)$ of the form $p(Z,\bar{Z})=0$, one automatically has $(pq)(Z,\bar{Z})=0$, for each polynomial $q$ with $\deg (pq)\leq n$.)  

Each \textit{singular} moment matrix $M(n)$ has at least one nontrivial
linear relation in its column space, and each such relation is naturally
associated with the zero set of a multivariable polynomial $p$. \ Consider now the ideal $%
\mathcal{I}$ generated by all the above mentioned polynomials. \ H.M. Mo\"eller \cite{Moe} and C. Scheiderer proved 
independently that $\mathcal{I}$ is a real radical ideal whenever $M(n)\geq 0$ (cf. \cite[Subsection 5.1, p. 203]{Lau3} and \cite{Lau2}). \
We also know that if $\mathcal{V}(\gamma )$ is finite, then $\mathcal{I}$ is
zero-dimensional, so results from algebraic geometry apply. 

We now recall a result that will allow us to convert a
given moment problem into a simpler, equivalent, moment problem. For $%
a,b,c\in \mathbb{C}$, $\left| b\right| \neq \left| c\right| $, let $\varphi
(z):=a+bz+c\bar{z}$ ($z\in \mathbb{C}$). Given $\gamma ^{(2n)}$, define $%
\tilde{\gamma}^{\left( 2n\right) }$ by $\tilde{\gamma}_{ij}:=\Lambda (\bar{\varphi}^{i}\varphi ^{j})$ ($0\leq i+j\leq 2n$), where $\Lambda$
denotes the Riesz functional associated with $\gamma $. It is
straightforward to verify that if $\Phi \left( z,\bar{z}\right) :=\left(
\varphi \left( z\right) ,\overline{\varphi \left( z\right) }\right) $, then $%
L_{\tilde{\gamma}}(p)=\Lambda \left( p\circ \Phi \right) $ for every $%
p\in \mathcal{P}_{n}$. (Note that for $p\left( z,\bar{z}\right) \equiv \sum
a_{ij}\bar{z}^{i}z^{j}$, $\left( p\circ \Phi \right) \left( z,\bar{z}\right)
=p\left( \varphi \left( z\right) ,\overline{\varphi \left( z\right) }\right)
\equiv \sum a_{ij}\overline{\varphi \left( z\right) }^{\,i}\varphi \left(
z\right) ^{j}$.)

\begin{proposition}
\label{lininv} (Invariance under degree-one transformations \cite[Proposition 1.7]{tcmp6}.) \ Let $%
M(n)$ and $\tilde{M}(n)$ be the moment matrices associated with $\gamma$ and 
$\tilde{\gamma}$, and let $J\hat{p}:=\widehat{p\circ\Phi} (p\in \mathcal{P}_{n})$. Then the following statements hold.

\begin{enumerate}
\item  \label{lininv(1)}$\tilde{M}(n)=J^{\ast}M(n)J$.

\item  \label{lininv(2)}$J$ is invertible.

\item  \label{lininv(3)}$\tilde{M}(n)\geq0\Leftrightarrow M(n)\geq0$.

\item  \label{lininv(4)}$\operatorname{rank}\tilde{M}(n)=\operatorname{rank}M(n)$.

\item  \label{lininv(5)}The formula $\mu=\tilde{\mu}\circ\Phi$ establishes a
one-to-one correspondence between the sets of representing measures for $%
\gamma$ and $\tilde{\gamma}$, which preserves measure class and cardinality
of the support; moreover, $\varphi(\operatorname{supp}\mu )=\operatorname{supp}\tilde{%
\mu}$.

\item  \label{lininv(6)}$M\left( n\right) $ admits a flat extension if and
only if $\tilde{M}\left( n\right) $ admits a flat extension.

\item  \label{lininv(7)}For $p\in\mathcal{P}_{n}$, $p\left( \tilde{Z},%
\Tilde{\Bar{Z}}%
\right) =J^{\ast}\left( \left( p\circ\Phi\right) \left( Z,\bar{Z}\right)
\right) $.
\end{enumerate}
\end{proposition}

\subsection{\textbf{Statement of the Main Result}}

In this paper we initiate the study cubic column relations associated with
finite algebraic varieties, that is, the case when  
$M(3)\geq 0$, $M(2)>0$, with $\operatorname{card} \mathcal{V}(M(3))<\infty$. \ (Cubic column relations with infinite variety have already appeared in \cite{Fianew}.) \ We begin by stating the general solution of the singular \textit{quartic} moment problem.

\begin{theorem}
\label{quartic1} (cf. \cite{tcmp2}, \cite{tcmp4}, \cite{tcmp7}, \cite{tcmp9}, 
\cite{FiaNie}) \ Let $p \in \mathbb{C}[z,\bar{z}]$, with $\deg p \leq 2$. \ Then $\gamma ^{(2n)}$ has a representing measure supported
in the curve $p(z,\bar{z})=0$ if and only if $M(n)$ has a column dependence
relation $p(Z,\bar{Z})=0$, $M(n) \geq 0$, $M(n)$ is recursively
generated, and $r\leq v$.
\end{theorem}

In view of Theorem \ref{quartic1}, we can
always assume that $M(2)$ is invertible. \ Since $M(3)$ is a square matrix of size $10$, and since we focus on the case of
finite algebraic variety, the possible values for $\operatorname{rank} \; M(3)$ are $7$
and $8$. \ (While $\operatorname{rank} \; M(3)=9$ is theoretically possible, the associated algebraic
variety, which must have at least $9$ points to conform with the inequality $r \leq v$, would be infinite, being 
determined by a single column dependence relation.) \ When $r=v=7$, we focus on the case of 
a column relation given by an harmonic polynomial $q(z,\bar{z}):=f(z)-%
\overline{g(z)}$, where $f$ and $g$ are analytic polynomials, and $\deg \
q=3 $. \ Using degree-one transformations and symmetric properties of such
polynomials, we reduce TMP to the case $Z^{3}=itZ+u\bar{Z}$, with $%
t,u\in \mathbb{R}$. \ Wilmhurst \cite{Wil}, Crofoot-Sarason \cite{CrSa} and
Khavinson-Swiatek \cite{KhSw} proved that for $\deg \ f=3$, we have $\operatorname{%
card}\ \mathcal{Z}(f(z)-\bar{z})\leq 7$. \ It immediately follows that a TMP
with such a cubic column relation can have at most $7$
points in its algebraic variety. \ 

We present below our main result (Theorem \ref{thmcubic}). \ First, we need some notation. \ For $u,t \in \mathbb{R}$, let 
\begin{equation}
q_{7}(z,\bar{z}):=z^{3}-itz-u\bar{z} \label{q7def}
\end{equation}
and
\begin{equation}
q_{LC}(z,\bar{z}):=i(z-i \bar{z})(\bar{z}z-u). \label{qLCdef}
\end{equation}

For $u,t \in \mathbb{R}$ and $q_{7}$ defined as in (\ref{q7def}), assume that $(u,t)$ is in the open cone $(0<\left\vert
u\right\vert <t<2\left\vert u\right\vert )$. \ Then $\operatorname{card}\ \mathcal{Z}%
(q_{7})=7$, as can be verified using \textit{Mathematica} \cite{Wol} (cf. Lemma \ref{sec3main}). \ In fact, this $7$-point set consists of the origin, two points equidistant
from the origin, located on the bisector $z=i\bar{z}$, and four points on a
circle, symmetrically located with respect to the bisector (cf. Figure \ref%
{figurecubic} and Lemma \ref{sec3main}). \ Moreover, there is a
cubic polynomial whose zero set consists precisely of the union of the
bisector and the circle, given by $q_{LC}(z,\bar{z}):=i(z-i \bar{z}%
)(\bar{z}z-u)$. \ The fact that $\mathcal{Z}(q_7) \subseteq \mathcal{Z}%
(q_{LC})$ is crucial.

\setlength{\unitlength}{1mm} \psset{unit=15mm} 

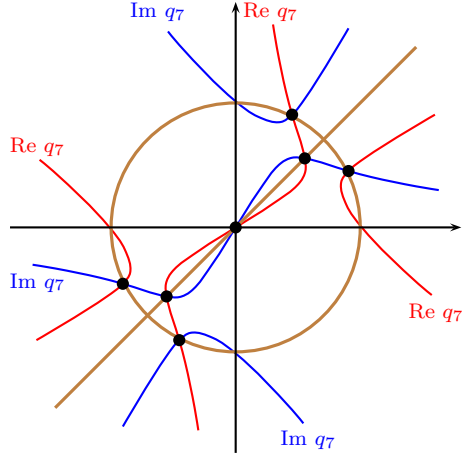
\begin{figure}[th]
\begin{center}

\begin{picture}(15,30)
%Re(q_7)
\pscurve[linecolor=red](-.331,-1.8)(-.364,-1.6)(-.447,-1.2)(-.5,-1)(-.612,-.612)(-.618,-.5)(-.569,-.4)(-.315,-.2)(0,0)(.315,.2)(.569,.4)(.618,.5)(.612,.612)(.5,1)(.447,1.2)(.364,1.6)(.331,1.8)

\pscurve[linecolor=red](1.738,-.6)(1.618,-.5)(1.118,0)(.981,.2)(.935,.4)(1,.5)(1.14,.612)(1.766,1)

\pscurve[linecolor=red](-1.766,-1)(-1.14,-.612)(-1,-.5)(-.935,-.4)(-.981,-.2)(-1.118,0)(-1.618,.5)(-1.738,.6)

%Im(q_7)

\pscurve[linecolor=blue](-1.8,-.331)(-1.6,-.364)(-1.2,-.447)(-1,-.5)(-.612,-.612)(-.5,-.618)(-.4,-.569)(-.2,-.315)(0,0)(.2,.315)(.4,.569)(.5,.618)(.612,.612)(1,.5)(1.2,.447)(1.6,.364)(1.8,.331)
\pscurve[linecolor=blue](-.6,1.738)(-.5,1.618)(0,1.118)(.2,.981)(.4,.935)(.5,1)(.612,1.14)(1,1.766)
\pscurve[linecolor=blue](-1,-1.766)(-.612,-1.14)(-.5,-1)(-.4,-.935)(-.2,-.981)(0,-1.118)(.5,-1.618)(.6,-1.738)

%qlc

\pscircle[linecolor=brown, linewidth=1.2pt](0,0){1.118}
\psline[linecolor=brown, linewidth=1.2pt](-1.6,-1.6)(1.6,1.6)

\psset{linecolor=black}
\qdisk(1,0.5){2.3pt}
\qdisk(0.5,1){2.3pt}
\qdisk(-1,-0.5){2.3pt}
\qdisk(-0.5,-1){2.3pt}
\qdisk(0.612372,0.612372){2.3pt}
\qdisk(-0.612372,-0.612372){2.3pt}
\qdisk(0,0){2.3pt}

\put(-30,10){\color{red}{\footnotesize Re $q_7$}}
\put(-14,28){\color{blue}{\footnotesize Im $q_7$}}
\put(23,-12){\color{red}{\footnotesize Re $q_7$}}
\put(-30,-8){\color{blue}{\footnotesize Im $q_7$}}
\put(1,28){\color{red}{\footnotesize Re $q_7$}}
\put(6,-29){\color{blue}{\footnotesize Im $q_7$}}

%axis
\psline{->}(-2,0)(2,0)
\psline{->}(0,-2)(0,2)

\end{picture}

\end{center}

\vspace{75pt}

\caption{{\footnotesize The $7$-point set $\mathcal{Z}(q_{7})$. \ (The {\color{brown} circle} has radius $\sqrt{u}$.)}}
\label{figurecubic}
\end{figure}

\begin{theorem}
\label{thmcubic} \ Let $M(3)\geq 0\,$, with $M(2)>0$ and $q_{7}(Z,\bar{Z})=0$. \ For $u,t \in \mathbb{R}$ and $q_{7}, q_{LC}$ defined as in (\ref{q7def}) and (\ref{qLCdef}), assume that $(0<\left\vert
u\right\vert <t<2\left\vert u\right\vert )$. \ The following statements are equivalent. \newline
(i) \ There exists a representing measure for $M(3)$. \newline
(ii) \ 
\begin{equation*}
\left\{ 
\begin{array}{ccc}
\Lambda (q_{LC}) & = & 0 \phantom{.} \\ 
\Lambda (zq_{LC}) & = & 0. %\label{condition}%
\end{array}%
\right.
\end{equation*}
(iii)  
\begin{equation*}
\left\{ 
\begin{array}{ccc}
\emph{Re}\; \gamma_{12}-\emph{Im}\; \gamma_{12}&=&u(\emph{Re}\; \gamma_{01}-\emph{Im}\; \gamma_{01}) \\ 
\gamma_{22}&=&(t+u)\gamma_{11}-2u \; \emph{Im}\; \gamma_{02} \nonumber%
\end{array}%
\right.
\end{equation*}
(iv) \ $q_{LC}(Z,\bar{Z})\equiv \bar{Z}^2 Z+i\bar{Z}Z^2-u\bar{Z}-iuZ=0$.
\end{theorem}

Since we are dealing with an extremal moment problem, in order to prove Theorem \ref%
{thmcubic} we need to verify Consistency (see \ref{C2}), which we do with the help of Lemma \ref{rep}.

\medskip \textit{Acknowledgments}. \ The authors are indebted to Professors Lawrence Fialkow and Frank Sottile for several discussions related 
to the topics in this paper. \ The authors are also deeply grateful to the referee for many suggestions that led to significant improvements in the presentation. \ Most of the examples, and some of the
proofs in this paper, were obtained with the aid of the software tool 
\textit{Mathematica} \cite{Wol}.

\section{Cubic Column Relations} \label{Section3}

Since we know how to solve the singular quartic moment problem \cite{tcmp6}, without loss of generality we will hereafter assume $M(2)>0$. \ We first recall a result from \cite{tcmp2}.

\begin{theorem} \label{quarticthm}
(\cite[Theorem 3.1]{tcmp2}) \ If $M(n)$ admits a column relation of the form $Z^{k}=p_{k-1}(Z,%
\bar{Z})\;\;$($1\leq k\leq \left[ \frac{n}{2}\right] +1$ and $\deg
\;p_{k-1}\leq k-1$), then $M(n)$ admits a flat extension $M(n+1)$, and
therefore a representing measure.
\end{theorem}

Now, if $k=3$, Theorem \ref{quarticthm} can be used only if 
$n\geq 4$. \ Thus, one strategy
is to somehow extend $M(3)$ to $M(4)$ and preserve 
the column relation $Z^{3}=p_{2}(Z,\bar{Z})$. 
\ This requires first checking that the $B$ block can be written as $M(3)W$ for some $W$, and then verifying that the $C$ block in the
extension satisfies the Toeplitz condition, something highly nontrivial. \ (A concrete attempt using \textit{Mathematica} \cite{Wol} led to memory overflow.) 

On the other hand, since we always assume that $M(2)$ is positive and invertible, and since the flat extension case ($\operatorname{rank} \; M(3) = \operatorname{rank} \; M(2)$) is well known, the first nontrivial case of $M(3)$ with finite variety arises when $\operatorname{rank} \; M(3)=7$. \ Now, since a soluble TMP requires $r \le v$, the algebraic variety of a soluble TMP needs to have a minimum of $7$ points. \ In other words, when $r=7$, we either have $v \ge 7$ or no representing measure. \ Now, given a cubic polynomial $p(z) \equiv z^3+bz^2+cz+d$, the substitution $w=z+b/3$ (which produces an equivalent TMP by Proposition \ref{lininv}) transforms it into $q(w) \equiv w^3+\tilde{c}w+\tilde{d}$; thus, without loss of generality, we always assume that there is no quadratic term in the analytic piece.    

Based on the previous considerations, we would like to focus our study on the case of harmonic polynomials, that is, polynomials of the form $q(z,\bar{z}):=f(z)-\overline{g(z)}$, with $\deg \;q=3$. \ In the case when $g(z)\equiv z$, we have

\begin{lemma} \label{sarason}
(\cite{Wil}, \cite{CrSa}, \cite{KhSw}) \ 
\[
\operatorname{card}\;\mathcal{Z}(f(z)-\bar{z})\leq 7.
\]
\end{lemma}

Observe that B\'ezout's Theorem predicts $\operatorname{card}\;\mathcal{Z}(f(z)-\bar{z})\leq 9$, so Lemma \ref{sarason} produces a better upper bound for the number of zeros. \ However, to get at least $7$ points is not generally easy, because most complex cubic harmonic polynomials have $5$ or fewer zeros. \ One way to maximize the number of zeros is to impose 
symmetry conditions on the zero set $K$. \ Also, for a polynomial of the form $z^{3}+\alpha z+\beta \bar{z}$, it is clear that $0 \in K$ and that $-z \in K$ whenever $z \in K$. \ Another natural condition is to require that $K$ be symmetric with respect to the line $y=x$, which in complex notation is $z=i \bar{z}$. \ When this is required, we obtain $\alpha \in i \mathbb{R}$ and $\beta \in \mathbb{R}$. \ Thus, the column relation becomes $Z^3=itZ+u \bar{Z}$, with $t,u \in \mathbb{R}$. 
 
Under these conditions, one needs to find only two points, one on the line $y=x$, the other outside that line. \ We thus consider the harmonic polynomial $q_{7}(z,\bar{z}):=z^{3}-itz-u\bar{z}$, with $u,t \in \mathbb{R}$. 
  
\begin{lemma} \label{sec3main}
Let $q_7$ be as above, and assume $0<\left| u\right| <t<2\left| u\right| $. \ Then $\operatorname{card}\;\mathcal{Z}(q_{7})=7$. \ 
In fact, 
\begin{equation} \label{expression}
\mathcal{Z}(q_{7})=\{0,\pm(p+iq),\pm(q+ip),\pm(r+ir) \},
\end{equation}
where $p,q,r>0$, $t=4pq$, $u=p^{2}+q^{2}$ and $r^{2}=\frac{t-u}{2}$ (cf. Figure \ref{picture4}).  
\end{lemma}

{\bf Proof.} 
We begin with a simple observation: for any pair of positive numbers $u$ and $t$ such that $0<u<t<2u$ one can always find a unique pair of positive numbers $p$ and $q$ such that $u=p^2+q^2$ and $t=4pq$. \ For, consider the functions $P(\theta):=\sqrt{u} \cos \theta$ and $Q(\theta):=\sqrt{u} \sin \theta$, on the interval $(\frac{\pi}{12},\frac{\pi}{4})$. \ It is straightforward to verify that $P(\theta)^2+Q(\theta)^2=u$ and that $4P(\theta)Q(\theta)=2u \sin(2 \theta)$, so that $4P(\theta)Q(\theta)$ maps $(\frac{\pi}{12},\frac{\pi}{4})$ onto $(u,2u)$.  It follows that any positive number $t$ between $u$ and $2u$ can be uniquely represented as $4pq$, with $p^2+q^2=u$. \ Also, observe that neither $p$ nor $q$ can be zero, and that $p \ne q$. 

Next, we identify the two real polynomials Re $q_7 =x^3-3xy^2+ty-ux$ and 
Im $q_7=-y^3+3x^2y-tx+uy$ (whose graphs are shown in Figure \ref{figurecubic}),  and calculate $Res(x):=\text{Resultant}(\text{Re~} q_7, \text{Im~} q_7,y)$, which is the determinant of the associated Sylvester matrix, i.e.,
\begin{eqnarray}
Res(x)&=&\det \left( 
\begin{array}{ccccc}
-3x & t & x^{3}-ux & 0 & 0 \\ 
0 & -3x & t & x^{3}-ux & 0 \\ 
0 & 0 & -3x & t & x^{3}-ux \\ 
-1 & 0 & 3x^{2}+u & -tx & 0 \\ 
0 & -1 & 0 & 3x^{2}+u & -tx%
\end{array}%
\right) \nonumber \\
&=& x\left( 2x^{2}+u-t \right) \left( 2x^{2}+u+t \right) \left(
16x^{4}-16ux^{2}+t^{2}\right) \label{resultant}.
\end{eqnarray}
As shown in real algebraic geometry, the resultant detects the common zeros of Re $q_7$ and Im $q_7$ \cite{CLO1}. \  From (\ref{resultant}), we immediately observe that:
\begin{itemize}
\item[(1)] One zero of $q_7$ (the origin) comes from the linear factor $x$;
\item[(2)] Two zeros of $q_7$ come from the factors $ 2x^{2}+u-t$ and $ 2x^{2}+u+t$ (which obviously cannot be simultaneously zero);
\item[(3)] Four zeros (which are necessarily located outside the bisector $z=i\bar{z}$) come from the factor $16x^{4}-16x^{2}u+t^{2}$.
\end{itemize}
It is then clearly sufficient to prove the result for the case $u>0$ which, using (2) above, consists of analyzing the factor $ 2x^{2}+u-t$. \ Thus, the condition  $u<t$ is essential to have two points on the bisector $z=i\bar{z}$ (besides the origin). \ It remains to investigate (3) above. \ Toward this end, consider $16x^4-16x^{2}u+t^{2}=0$. \ Then
\begin{eqnarray*}
x^2 = \frac{2u\pm \sqrt{4u^2-t^2}}{4},
\end{eqnarray*}
where the right hand side is always positive under the second necessary condition $4u^2-t^2>0$. \ Now recall that there exists a unique pair $(p,q)$ of positive, distinct numbers, such that $u=p^2+q^2$ and $t=4pq$. \ Thus, the expression $4u^2-t^2$ equals $4(p^2+q^2)-16p^2q^2=4(p^2-q^2)^2$. \ It now follows easily that $x=p$ or $x=q$, depending on whether $p>q$ or $p<q$. \ With this at hand, it is straightforward to identify the $7$ points, as listed in (\ref{expression}) (see Figure \ref{picture4}). \qed

\begin{remark}
The fact that $q_{7}$ has the maximum number of zeros 
predicted by Lemma \ref{sarason}
is significant to us, in that each sextic TMP with 
invertible $M(2)$ and a
column relation of the form $q_{7}(Z,\bar{Z})=0$ either 
does not admit a
representing measure or is necessarily extremal.
\end{remark}
 
As a consequence, the existence of a representing measure will be
established once we prove that such a TMP is consistent. \ 
This means that for each polynomial $p$ of degree at most $6$ 
that vanishes on $\mathcal{Z}(q_{7})$ we must verify that $\Lambda(p)=0$.

\subsection{\textbf{The Hidden Column Relation}}

Since $\operatorname{rank} M(3)=7$, there must be another column relation 
besides $q_7(Z,\bar{Z})=0$ (and its conjugate). \ Clearly the columns
\begin{eqnarray*}
1,Z,\bar{Z},Z^2,\bar{Z}Z,\bar{Z}^2,\bar{Z}Z^2
\end{eqnarray*}
must be linearly independent (otherwise $M(3)$ would be a flat extension of $M(2)$), so the 
new column relation must involve both $\bar{Z}Z^2$ and $\bar{Z}^2Z$. \ In what follows, we prove that the ``hidden column relation" $R(Z,\bar{Z})=0$ is uniquely determined by the zero set of $q_7$. \ With the notation of Theorem \ref{sec3main}, let $P_0:=0$, $P_1:=p+iq$, $P_2:=-P_1$, $P_3:=q+ip$, $P_4:=-P_3$, $P_5:=r+ir$ and $P_6:=-P_5$. \ Let $R$ denote the polynomial giving rise to the column relation $R(Z,\bar{Z})=0$. \ Since the coefficients of $\bar{z} z^2$ and $\bar{z}^2 z$ are nonzero, without loss of generality we can assume that the coefficient of $\bar{z}^2 z$ is $1$. \ That is,
$$R(z,\bar{z}) \equiv \bar{z}^2 z +a_{12}\bar{z}z^2+a_{20}\bar{z}^2+a_{11}\bar{z}z+a_{02}z^2+a_{10}\bar{z}+a_{01}z+a_{00}.$$ \ 
By evaluating at $P_0$, it is easy to see that $a_{00}=0$. \ Moreover, the evaluations at $P_i (i=1,\cdots,6)$ can be presented as 
\begin{equation}
\label{vander6}
\left(
\begin{array}{c}
R(P_1) \\
R(P_2) \\
R(P_3) \\
R(P_4) \\
R(P_5) \\
R(P_6) 
\end{array}
\right)
\equiv
\left(
\begin{array}{c}
\bar{P_1}^2 P_1 \\
\bar{P_2}^2 P_2 \\
\bar{P_3}^2 P_3 \\
\bar{P_4}^2 P_4 \\
\bar{P_5}^2 P_5 \\
\bar{P_6}^2 P_6 
\end{array}
\right)
+
\left(
\begin{array}{cccccc}
P_1 & \bar{P_1} & P_1^2 & \bar{P_1} P_1 & \bar{P_1}^2 & \bar{P_1} P_1^2 \\
P_2 & \bar{P_2} & P_2^2 & \bar{P_2} P_2 & \bar{P_2}^2 & \bar{P_2} P_2^2 \\
P_3 & \bar{P_3} & P_3^2 & \bar{P_3} P_3 & \bar{P_3}^2 & \bar{P_3} P_3^2 \\
P_4 & \bar{P_4} & P_4^2 & \bar{P_4} P_4 & \bar{P_4}^2 & \bar{P_4} P_4^2 \\
P_5 & \bar{P_5} & P_5^2 & \bar{P_5} P_5 & \bar{P_5}^2 & \bar{P_5} P_5^2 \\
P_6 & \bar{P_6} & P_6^2 & \bar{P_6} P_6 & \bar{P_6}^2 & \bar{P_6} P_6^2 
\end{array}
\right)
\left(
\begin{array}{c}
a_{01} \\
a_{10} \\
a_{02} \\
a_{11} \\
a_{20} \\
a_{12} 
\end{array}
\right)
=
\left(
\begin{array}{c}
0 \\
0 \\
0 \\
0 \\
0 \\
0 
\end{array}
\right) .
\end{equation}
 
A calculation using \textit{Mathematica} \cite{Wol} shows that the determinant of the $6 \times 6$ matrix in (\ref{vander6}) is $128(1+i)(p-q)^4(p+q)^2r^2(p^2+q^2-2r^2)$. \ Since $p^2+q^2=u$ and $2r^2=t-u$, we see that the last factor is $2u-t$, and as a result the above mentioned determinant is different from zero. \ It follows that there exists exactly one monic polynomial $R$ vanishing in the $7$-point set $\mathcal{Z}(q_7)$. \ On the other hand, it is not hard to see that the polynomial 
\begin{equation} \label{qLC2}
q_{LC}(z,\bar{z}):=i(z-i \bar{z})(\bar{z}z-u) \; \; (\textrm{cf.} \; (\ref{qLCdef}))
\end{equation}
vanishes in the zero set of $q_7$, and it is monic, so it must be $R$. \ We have thus found the ``hidden column relation": it is
$$
q_{LC}(Z,\bar{Z}) \equiv \bar{Z}^2 Z+i\bar{Z}Z^2-u\bar{Z}-iuZ=0.
$$

\begin{remark}
Since  
\begin{eqnarray*}
q_{LC}(z,\bar{z})&=&i(z-i\bar{z})(\bar{z}z-u), \\ 
\end{eqnarray*}
it is straightforward to see that the zero set of $q_{LC}$ is the union of a line and a 
circle, and that $\mathcal{Z}(q_7) \subseteq \mathcal{Z}(q_{LC})$ (see Figure \ref{picture4}).
\end{remark}

\setlength{\unitlength}{1mm} \psset{unit=1mm}
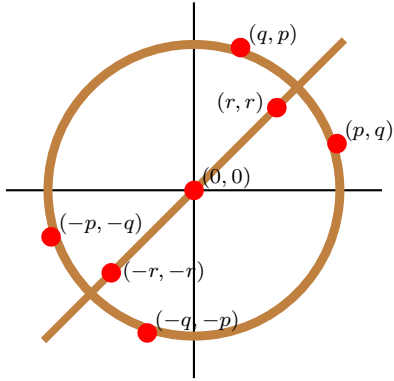
\begin{figure}[th]
\vspace{48pt}
\begin{center}
\begin{picture}(60,30)

\psline(5,20)(55,20)
\psline(30,-5)(30,45)
\psline[linecolor=brown,linewidth=.1cm](10,-0)(50,40)
\pscircle[linecolor=brown,linewidth=.12cm](30,20){20}

\put(30,20){\pscircle*[linecolor=red](0,0){1.3}}
\put(41,31){\pscircle*[linecolor=red](0,0){1.3}}
\put(19,9){\pscircle*[linecolor=red](0,0){1.3}}
\put(49,26.2){\pscircle*[linecolor=red](0,0){1.3}}
\put(36.2,39){\pscircle*[linecolor=red](0,0){1.3}}
\put(11,13.8){\pscircle*[linecolor=red](0,0){1.3}}
\put(23.8,1){\pscircle*[linecolor=red](0,0){1.3}}

\put(31,21){\footnotesize{$(0,0)$}}
\put(33,31){\footnotesize{$(r,r)$}}
\put(20.5,8.5){\footnotesize{$(-r,-r)$}}
\put(50,27.2){\footnotesize{$(p,q)$}}
\put(37.2,40){\footnotesize{$(q,p)$}}
\put(12,14.8){\footnotesize{$(-p,-q)$}}
\put(24.8,2){\footnotesize{$(-q,-p)$}}

\end{picture}
\end{center}

\caption{The sets ${\color{red} \mathcal{Z}(q_{7})}$ and 
{\color{brown} $\mathcal{Z}(q_{LC})$}; here $r=\sqrt{\frac{t-u}{2}}$, 
$p=\frac{1}{2}(2u+\sqrt{4u^2-t^2})$ and $p^2+q^2=u$.}
\label{picture4}
\end{figure}

%%%%%%%%%%%%%%%%%

\section{Main Theorem} \label{proofs}

We are ready to prove Theorem \ref{thmcubic}, which we restate for the reader's convenience.

\begin{theorem}
\label{thmcubic2} \ Let $M(3)\geq 0\,$, with $M(2)>0$ and $q_{7}(Z,\bar{Z})=0$. \ For $u,t \in \mathbb{R}$ and $q_{7}, q_{LC}$ defined as in (\ref{q7def}) and (\ref{qLCdef}), assume that $0<\left\vert
u\right\vert <t<2\left\vert u\right\vert $. \ The following statements are equivalent. \newline
(i) \ There exists a representing measure for $M(3)$. \newline
(ii) \ 
\begin{equation*}
\left\{ 
\begin{array}{ccc}
\Lambda (q_{LC}) & = & 0 \phantom{.} \\ 
\Lambda (zq_{LC}) & = & 0. \label{condition}%
\end{array}%
\right.
\end{equation*}
(iii)  
\begin{equation*}
\left\{ 
\begin{array}{ccc}
\emph{Re}\; \gamma_{12}-\emph{Im}\; \gamma_{12}&=&u(\emph{Re}\; \gamma_{01}-\emph{Im}\; \gamma_{01}) \\ 
\gamma_{22}&=&(t+u)\gamma_{11}-2u \; \emph{Im}\; \gamma_{02} \nonumber%
\end{array}%
\right.
\end{equation*}
(iv) \ $q_{LC}(Z,\bar{Z})=0$.
\end{theorem}

To prove Theorem \ref{thmcubic2}, we will need the following auxiliary result.

\begin{lemma} (Representation of Polynomials) \label{rep}
\ For $u$ and $t$ as in Theorem \ref{thmcubic2}, let $\mathcal{P}_{6}:=\{p\in \mathbb{C}_{6}[z,\bar{z}]:p|_{\mathcal{Z}%
(q_{7})}\equiv 0\}$   and let $\mathcal{I}:=\{p\in \mathbb{C}_{6}[z,\bar{z}%
]:p=fq_{7}+g\bar{q}_{7}+hq_{LC}$ for some $f,g,h\in \mathbb{C}_{3}[z,\bar{z}%
]\}$.   \ Then $\mathcal{P}_{6}=\mathcal{I}$. 
\end{lemma}

{\bf Proof.} \ 
Clearly, $\mathcal{I}\subseteq \mathcal{P}_{6}$. \ We shall show that $\dim
\;\mathcal{I}=\dim \;\mathcal{P}_{6}$.   \ Let $T:\mathbb{C}%
^{30}\longrightarrow \mathbb{C}_{6}[z,\bar{z}]$ be given by 
\begin{eqnarray*}
(a_{00},\cdots ,a_{30},b_{00},\cdots ,b_{30},c_{00},\cdots ,c_{30})
&\longmapsto &(a_{00}+a_{01}z+a_{10}\bar{z}+\cdots+a_{30}\bar{z}^{3})q_{7} \\
&&+(b_{00}+b_{01}z+b_{10}\bar{z}+\cdots+b_{30}\bar{z}^{3})\bar{q}_{7} \\
&&+(c_{00}+c_{01}z+c_{10}\bar{z}+ \cdots +c_{30}\bar{z}^{3})q_{LC}.
\end{eqnarray*}

Recall that $30 = \dim \; \mathbb{C}^{30}=\dim \; \ker \; T + \dim \; $ran$ \; T$ (by the Fundamental Theorem of Linear Algebra), and observe that $\mathcal{I}=$ ran $\; T$,   so that 
$\dim \;\mathcal{I}=\operatorname{%
rank}\;T$.

To find $\operatorname{rank}\;T$, we first determine $\dim
\;\ker \;T$. \ Using Gaussian elimination (with the aid of \textit{Mathematica} \cite{Wol}), we can prove that 
$\dim \;\ker \;T=9$ 
whenever $ut\neq 0$.   \ It follows that $\operatorname{rank}\;T=30-9=21$, 
that is, $\dim \;\mathcal{I}=21$. \

Now consider the evaluation map $S:\mathbb{C}_{6}[z,\bar{z}]\longrightarrow 
\mathbb{C}^{7}$ given by 
\begin{eqnarray*}
S(p(z,\bar{z}))&:=&(p(w_{0},\bar{w}_{0}),p(w_{1},\bar{w}_{1}),
p(w_{2},\bar{w}_{2}), \\
&&p(w_{3},\bar{w}_{3}),p(w_{4},\bar{w}_{4}),p(w_{5},\bar{w}_{5}),
p(w_{6},\bar{w}_{6})).
\end{eqnarray*}

We know that $\dim \; \ker \; S + \dim \; $ran $\; S=\dim \; \mathbb{C}_{6}[z,\bar{z}]=28$.
Using Lagrange Interpolation, we can verify that $S$ 
is onto, i.e., $\operatorname{rank} \; S=7$. \ For, if we define
$$
\ell_j(z,\bar z):=\prod_{\stackrel{i=0}{i\neq j}}^6 \frac{z-w_i}{w_j - w_i},
$$
then $S(\ell_j(z,\bar z))=e_{j+1}$, where $e_{j}$ is the Euclidean basis element in $\mathbb C^7$ for $j=0,\cdots, 6$. \  Thus, $\{e_1,\cdots,e_7\} \in \text{ran~} S$. \ 

Now, it is straightforward to note that $\ker \;S=\mathcal{P}_{6}$, and since $\dim \;\mathbb{C}_{6}[z,\bar{z%
}]=28$, it follows that $\dim \;\ker \;S=21$, and a fortiori that $\dim \;%
\mathcal{P}_{6}=21$. 

  \ Therefore, $\dim \; \mathcal{I} = 21 = \dim \; \mathcal{P}_{6}$, and since $\mathcal{I} \subseteq \mathcal{P}_{6}$, we have established that $\mathcal{I}=\mathcal{P}_{6}$, as desired. \qed

\medskip

{\bf Proof of Theorem \ref{thmcubic2}.} \ Conditions (ii), (iii) and (iv) are easily seen to be equivalent, so we focus on the proof of (i) $\Leftrightarrow$ (ii).

\noindent $(\Longrightarrow )$ \ Let $\mu $ be a representing measure. \ 
We know that $%
7\leq \operatorname{rank}\;M(3)\leq \operatorname{card}\;\operatorname{supp}\;\mu \leq 
\operatorname{card}\;\mathcal{Z}(q_{7})=7$, so that $\operatorname{supp}\;\mu =%
\mathcal{Z}(q_{7})$ and $\operatorname{rank}\; M(3)=7$. \ Thus,
\begin{equation*}
\Lambda(q_7)=\int q_7 \; d\mu=0.
\end{equation*}
Similarly, since 
$\operatorname{supp} \; \mu \subseteq \mathcal{Z}(q_{LC})$, we also have
\begin{equation*}
\Lambda(q_{LC})=\Lambda(zq_{LC})=0,
\end{equation*}
as desired.

\noindent $(\Longleftarrow )$ \ We will apply \cite[Theorem 4.2]{tcmp11}, that is, we will verify (\ref{C2}) (Consistency of $M(3)$). \ First, recall that $q_7(Z,\bar{Z})=0$ or, equivalently, $M(3)\hat{q_7}=0$, so that $\Lambda(fq_7)=\left\langle M(3)\hat{q_7},\hat{f}\right\rangle=0 \; (\textrm{for each polynomial~} f \; \textrm{in~} z \; \textrm{and~} \bar{z} \; \textrm{of degree at most ~} 3)$, from which it follows that the Riesz functional $\Lambda$ annihilates all polynomials of the form $fq_7$ if deg $f \le 3$. \ Using $q_7$ and (\ref{qLC2}), we can prove that 
$\Lambda (\bar{z}^{i}z^{j}q_{LC})=0$ for all $0\leq i+j\leq 3$. \ 
 
For example, 
\begin{eqnarray*}
\bar{z}q_{LC}-izq_{LC} &=&(\bar{z}-iz)(\bar{z}^{2}z+i\bar{z}z^{2}-iuz-u\bar{z%
})  \nonumber \\
&=&-uz^{2}+\bar{z}z^{3}-u\bar{z}^{2}+\bar{z}^{3}z  \nonumber \\
&=&-uz^{2}+\bar{z}(q_7+itz+u\bar{z})-u\bar{z}^{2}+(\bar{q_7}-it\bar{z}+uz)z \\
&=&\bar{z}q_7+\bar{q_7}z,  \nonumber
\end{eqnarray*}
and therefore $\Lambda (\bar{z}q_{LC})=i\Lambda (zq_{LC})=0$. \ Rather than testing all possible monomials, we proceed as follows. \  
First, we can reduce each polynomial $\bar z^i z^j q_{LC}$ for $2\leq i+j\leq3$ by using the identity $z^3=it z+u\bar z$ mod ker $\Lambda$. \ For example,
\begin{eqnarray*}
              z^{2}q_{LC}&=& z^{2}(\bar{z}^{2}z+i\bar{z}z^{2}-iuz-u\bar{z}) \\
    & =&\bar{z}^{2}z^{3}+i\bar{z}z^{4}-u\bar{z}z^{2}-iuz^{3}\\
    & =&\bar{z}^{2}(q_7+itz+u\bar{z})+i\bar{z}z(q_7+itz+u\bar{z})-u\bar{z}z^{2}-iu(q_7+itz+u\bar{z})\\
&= & i ( t+u) q_{LC}+hq_7,
           \end{eqnarray*}
for some polynomial $h$. \ (Indeed, $h(z,\bar{z}) \equiv \bar{z}^2 +i\bar{z}z-iu$.) \ Similarly, we have
\begin{eqnarray*}
 \bar z z q_{LC} &=&-(t+u)q_{LC} \textrm{~mod ker~} \Lambda;\\
 z^3 q_{LC} &=&( t+u) \bar{z}q_{7} \textrm{~mod ker~} \Lambda; \\
 \bar z z^2 q_{LC} &=& i( t+u) \bar{z}q_{7} \textrm{~mod ker~} \Lambda.
\end{eqnarray*}
Second, note that the Riesz functional is invariant under conjugation, and that $q_{LC}=\overline{q_{LC}}$. \ Now, we compute the Riesz functional for each polynomial of the form $\bar{z}^iz^jq_{LC}$, for $0 \le i+j \le 3$. \ Suppose that $\Lambda (q_{LC})=0$ and $\Lambda (z q_{LC})=0$. \ Then,
\begin{eqnarray*}
  \Lambda (\bar{z}q_{LC})&=&\Lambda (\bar{z}
\overline{q_{LC}})=\overline{\Lambda (zq_{LC})}=0;\\
\Lambda (z^2 q_{LC})&=&i(t+u)\Lambda (q_{LC})=0 ;\\
\Lambda (\bar z z q_{LC})&=&-i(t+u)\Lambda (q_{LC})=0;\\
\Lambda (\bar z^2 q_{LC})&=&\Lambda (\bar z^2 \overline{q_{LC}})= \overline{\Lambda (z^2 q_{LC})}=0;\\
 \Lambda (z^3 q_{LC})&=&(t+u) \Lambda (\bar z q_{LC})=0;\\
\Lambda (\bar z z^2 q_{LC})&=&i(t+u) \Lambda (\bar z  q_{LC})=0;\\
\Lambda (\bar z^2 z q_{LC})&=&\overline{ \Lambda (\bar z z^2  q_{LC})}=0;\\
\Lambda (\bar z^3 q_{LC})&=&\overline{\Lambda ( z^3 q_{LC})}=0.
\end{eqnarray*}

\ It follows
that for $f,g,h\in \mathbb{C}_{3}[z,\bar{z}]$ we have 
$\Lambda (fq_{7}+g\bar{q}_{7}+hq_{LC})=0$.   \ 
Consistency will be established once we show that all
degree-$6$ polynomials vanishing in $\mathcal{Z}(q_{7})$ are of the form 
$fq_{7}+g\bar{q}_{7}+hq_{LC}$. \ But this is the content of Lemma \ref{rep}. \qed

%%%%%%%%%%%%%%%%%%%%%%%%%%%%%%%%%%%%%%%%%%%%%%%%%%%%%%%%%%%%%%%%%%%%%%%%%%%%%%%%%%%%%%%%
%%%%%%%%%%%%%%%%%%%%%%%%%%%%%%%%%%%%%%%%%%%%%%%%%%%%%%%%%%%%%%%%%%%%%%%%%%%%%%%%%%%%%%%%

\section{\textbf{Some Concrete Examples}}

\begin{example}\label{q7-meas} 
(The case of a column relation associated with an irreducible cubic) \ Consider $M(3)(\gamma)$ given as
$$
\left(
\begin{array}{cccccccccc}
 1 & 0 & 0 & \frac{11 i}{14} & \frac{13}{14} & -\frac{11 i}{14} & 0 & 0 & 0 & 0 \\
 & & & & & & & & & \\
 0 & \frac{13}{14} & -\frac{11 i}{14} & 0 & 0 & 0 & \frac{7 i}{8} & \frac{59}{56} & -\frac{7 i}{8} & -\frac{23}{56} \\
 & & & & & & & & & \\
 0 & \frac{11 i}{14} & \frac{13}{14} & 0 & 0 & 0 & -\frac{23}{56} & \frac{7 i}{8} & \frac{59}{56} & -\frac{7 i}{8} \\
 & & & & & & & & & \\
 -\frac{11 i}{14} & 0 & 0 & \frac{59}{56} & -\frac{7 i}{8} & -\frac{23}{56} & 0 & 0 & 0 & 0  \\
 & & & & & & & & & \\
 \frac{13}{14} & 0 & 0 & \frac{7 i}{8} & \frac{59}{56} & -\frac{7 i}{8} & 0 & 0 & 0 & 0 \\
 & & & & & & & & & \\
 \frac{11 i}{14} & 0 & 0 & -\frac{23}{56} & \frac{7 i}{8} & \frac{59}{56} & 0 & 0 & 0 & 0 \\
 & & & & & & & & & \\
 0 & -\frac{7 i}{8} & -\frac{23}{56} & 0 & 0 & 0 & \frac{277}{224} & -\frac{227 i}{224} & -\frac{97}{224} & -\frac{61 i}{224} \\
& & & & & & & & & \\
 0 & \frac{59}{56} & -\frac{7 i}{8} & 0 & 0 & 0 & \frac{227 i}{224} & \frac{277}{224} & -\frac{227 i}{224} & -\frac{97}{224} \\
& & & & & & & & & \\
 0 & \frac{7 i}{8} & \frac{59}{56} & 0 & 0 & 0 & -\frac{97}{224} & \frac{227 i}{224} & \frac{277}{224} & -\frac{227 i}{224} \\
& & & & & & & & & \\
 0 & -\frac{23}{56} & \frac{7 i}{8} & 0 & 0 & 0 & \frac{61 i}{224} & -\frac{97}{224} & \frac{227 i}{224} & \frac{277}{224} \\
& & & & & & & & & 
\end{array}
\right).
$$
Using the Nested Determinants Test and Smul'jan's Theorem (Theorem \ref{smul}), it can be verified that $M(3)$ is positive semidefinite and is of rank $7$ with the three column relations 
$$Z^3= 2 i Z+\frac{5}{4} \bar Z,$$  
$$\bar Z^2 Z=i \frac{5}{4} Z + \frac{5}{4} \bar Z - i \bar Z Z^2,$$ 
and 
$$\bar Z^3= -2 i \bar Z+\frac{5}{4}  Z.$$ 
Notice that $t=2$ and $u=\frac{5}{4}$. \ The associated algebraic variety consists of exactly $7$ points: $z_1=0$, $z_2=1+\frac{1}{2}i$, $z_3=\frac{1}{2}+i$, $z_4=-1-\frac{1}{2}i$, $z_5=-\frac{1}{2}-i$, $z_6=\frac{\sqrt{6}}{4}+\frac{\sqrt{6}}{4}i$, and $z_7=-\frac{\sqrt{6}}{4}-\frac{\sqrt{6}}{4}i$. \ Thus $\gamma$ is extremal and the main theorem implies that $M(3)$ has a 7-atomic representing measure $\mu$. \ By the Flat Extension Theorem \cite[Remark 3.15, Theorem 5.4, Corollary 5.12, Theorem 5.13, and Corollary 5.15]{tcmp1}, the densities can be computed by solving a Vandermonde equation; the representing measure is given by $\mu=\sum_{i=1}^{7} \rho_i \delta_{z_i}$, where the densities $\rho_i=\frac{1}{7}$ for $1\leq i \leq 7$.
\end{example}

The next example shows that a column relation allowing 7 points does not guarantee the existence of a representing measure even though the moment matrix is positive semidefinite and recursively generated.

\begin{example}
(The case of a column relation with $7$ points and no representing measure) \ Consider $M(3)(\gamma)$ given as
$$\left(
\begin{array}{cccccccccc}
 1 & 0 & 0 & \frac{11 i}{14} & \frac{13}{14} & -\frac{11 i}{14} & 0 & 0 & 0 & 0 \\
 & & & & & & & & & \\
 0 & \frac{13}{14} & -\frac{11 i}{14} & 0 & 0 & 0 & \frac{7 i}{8} & \frac{21}{20} & -\frac{7 i}{8} & -\frac{23}{56} \\
 & & & & & & & & & \\
 0 & \frac{11 i}{14} & \frac{13}{14} & 0 & 0 & 0 & -\frac{23}{56} & \frac{7 i}{8} & \frac{21}{20} & -\frac{7 i}{8} \\
 & & & & & & & & & \\
 -\frac{11 i}{14} & 0 & 0 & \frac{21}{20} & -\frac{7 i}{8} & -\frac{23}{56} & 0 & 0 & 0 & 0 \\
 & & & & & & & & & \\
 \frac{13}{14} & 0 & 0 & \frac{7 i}{8} & \frac{21}{20} & -\frac{7 i}{8} & 0 & 0 & 0 & 0 \\
 & & & & & & & & & \\
 \frac{11 i}{14} & 0 & 0 & -\frac{23}{56} & \frac{7 i}{8} & \frac{21}{20} & 0 & 0 & 0 & 0 \\
 & & & & & & & & & \\
 0 & -\frac{7 i}{8} & -\frac{23}{56} & 0 & 0 & 0 & \frac{277}{224} & -\frac{161 i}{160} & -\frac{7}{16} & -\frac{61 i}{224} \\
 & & & & & & & & & \\
 0 & \frac{21}{20} & -\frac{7 i}{8} & 0 & 0 & 0 & \frac{161 i}{160} & \frac{277}{224} & -\frac{161 i}{160} & -\frac{7}{16} \\
 & & & & & & & & & \\
 0 & \frac{7 i}{8} & \frac{21}{20} & 0 & 0 & 0 & -\frac{7}{16} & \frac{161 i}{160} & \frac{277}{224} & -\frac{161 i}{160} \\
 & & & & & & & & & \\
 0 & -\frac{23}{56} & \frac{7 i}{8} & 0 & 0 & 0 & \frac{61 i}{224} & -\frac{7}{16} & \frac{161 i}{160} & \frac{277}{224} \\
 & & & & & & & & & 
\end{array}
\right)$$
Applying Smul'jan's Theorem (Theorem \ref{smul}), we know  $M(3)$ is positive semidefinite  if and only if $M(3)_{\mathcal B}$, the compression of $M(3)$ to the rows and columns indexed by the basis $\mathcal B$ for $\mathcal C_{M(3)}$, is positive semidefinite. \ Using \it{Mathematica} \cite{Wol}, we can verify that all nested determinants of $M(3)_{\mathcal B}$ are positive, and therefore it follows that  $M(3)$ is positive semidefinite. \ Row reduction via \emph{Mathematica} shows $M(3)$ has only two column relations ($Z^3= 2 i Z+\frac{5}{4} \bar Z$ and its conjugate $\bar{Z}^3= -2 i Z+\frac{5}{4} Z$) and hence, $M(3)$ has rank 8. \   As seen in Example \ref{q7-meas}, the zero set of the polynomial $z^3- 2 i z-\frac{5}{4} \bar z$ consists of $7$ points, and therefore the algebraic variety has at most 7 points. \ As a result, the variety condition (\ref{C3}) fails, and therefore there is no representing measure. \

\end{example}

We end this section by introducing another class of cubic harmonic polynomials with real coefficients only. \ We discovered this class independently of the previous class, using symmetry. \ Later, we learned that the two classes are indeed equivalent (at least from the perspective of TMP), via a degree-one transformation. \

\begin{corollary}\label{lem-another-q7}
Suppose ${M}(3)(\widehat{\gamma})$, the associated moment matrix of a moment sequence $\widehat{\gamma}$, is positive semidefinite and satisfies the column relation 
\begin{equation} \label{w}
W^{3}=2\alpha W-\beta \bar{W}
\end{equation}
for $0<\alpha<\left|\beta\right|<2\alpha$ and ${M}(2)(\widehat{\gamma})>0$. \ Then $\widehat{\gamma }$ has a representing measure if and only if
\begin{center}
$\Lambda(\hat q_{LC})=0$ and $ \Lambda(w \hat q_{LC})=0$,
\end{center}
where $\hat q_{LC}(w,\bar w):=\bar{w}^2 w-\bar{w} w^2+\beta w-\beta \bar{w}$.

%$\mathrm{\rm{Im}}(\widehat{\gamma }_{12})=\beta \mathrm{\rm{Im}}(%
%\widehat{\gamma }_{12})$ and  $\widehat{\gamma }_{22}=(2\alpha +\beta )%
%\widehat{\gamma }_{11}-2\beta \rm{Re}(\widehat{\gamma }_{20}).$
\end{corollary}

{\bf Proof.} \ We will prove this using the equivalence of TMP under
degree-one transformations. \ Consider the degree-one transformation %\linebreak
$$
w \equiv \varphi (z,\bar{z}):=(1+i)\bar{z}
$$
and let $M(3)$ be the moment matrix at the $(z,\bar{z})$ level. \ Now, we can transform a column relation in ${M}(3)(\widehat{\gamma})$ using $\varphi$, that is, $p(W,\bar{W})=J^{\ast }(p\circ \Phi )(Z,\bar{Z})$ as in Proposition \ref{lininv}. \ In the specific case of the column relation (\ref{w}), we easily obtain
$$
\bar{Z}^3=-i \alpha \bar{Z}+\frac{\beta}{2}Z
$$
(in the column space of $M(3)$) and therefore 
$$
Z^3=i\alpha Z+\frac{\beta}{2} \bar{Z}.
$$
We then take $t=\alpha$ and $u=\frac{\beta}{2}$ in Theorem \ref{thmcubic2} to get $0<\alpha<\left|\beta\right|<2\alpha$. \ Moreover, in this case the hidden column relation is 
$$
\bar{Z}^2 Z+i \bar{Z} Z^2 -u \bar{Z} - iuZ=0.
$$
It follows that, at the $(w,\bar{w})$ level, the hidden column relation is
$$
\bar{W}^2W-\bar{W}W^2-2u\bar{W}+2uW=0,
$$
i.e., $\hat q_{LC}(W,\bar W)=0$. \ This completes the proof. \qed

\medskip
The polynomial $r(w,\bar{w}):=w^{3}-2\alpha w+\beta \bar{w}$, which defines the column relation in (\ref{w}), has already appeared in the work of D. Khavinson and G. Neumann \cite{KhNe}.

\end{document}